\documentclass[12pt]{amsart}

\usepackage{amssymb}
\usepackage{latexsym}
\usepackage{amsmath}

\begin{document}

\def\c{{\mathbb C}}
\def\r{{\mathbb R}}
\def\q{{\mathbb Q}}
\def\z{{\mathbb Z}}
\def\i{{\mathbb I}}

\newtheorem{theo}{Theorem}[section]
\newtheorem{prop}{Proposition}[section]
\newtheorem{lemm}{Lemma}[section]
\newtheorem{defi}{Definition}[section]
\newtheorem{coro}{Corollary}[section]
\newtheorem{exam}{Example}[section]
\newtheorem{rema}{Remark}[section]
\newtheorem{conj}{Conjecture}[section]

\title{Free product formulae for quantum permutation groups}
\author{Teodor Banica}
\address{Departement de Mathematiques, Universite Paul Sabatier, 118 route de Narbonne, 31062 Toulouse, France}
\email{banica@picard.ups-tlse.fr}
\author{Julien Bichon}
\address{Laboratoire de Mathematiques Appliquees, Universite de Pau et des Pays de l'Adour, IPRA, Avenue de l'universite, 64000 Pau, France}
\email{bichon@univ-pau.fr}

\begin{abstract}
Associated to a finite graph $X$ is its quantum automorphism group $G(X)$. We prove a formula of type $G(X*Y)=G(X)*_wG(Y)$, where $*_w$ is a free wreath product. Then we discuss representation theory of free wreath products, with the conjectural formula $\mu (G*_wH)=\mu (G)\boxtimes\mu (H)$, where $\mu$ is the associated spectral measure. This is verified in two situations: one using free probability techniques, the other one using planar algebras.  
\end{abstract}
\maketitle

\section*{Introduction}

A quantum group is an abstract object, dual to a Hopf algebra. Finite quantum groups are those which are dual to finite dimensional Hopf algebras.

A surprising fact, first noticed by Wang in \cite{wa2}, is that the quantum group corresponding to the Hopf algebra $\c^*(\z_2*\z_2)$ has a faithful action on the set $\{1,2,3,4\}$. This quantum group, which is of course not finite, is a so-called quantum permutation group.

In general, a quantum permutation group $G$ is described by a special type of Hopf $\c^*$-algebra $A$, according to the heuristic formula $A=\c (G)$. See \cite{gr}, \cite{wa2}.

The simplest case is when $A$ is commutative. Here $G$ is a subgroup of the symmetric group $S_n$. This situation is studied by using finite group techniques.

In general $A$ is not commutative, and infinite dimensional. In this case $G$ is a non-classical, non-finite compact quantum group. There is no analogue of a Lie algebra in this situation, but several representation theory methods, due to Woronowicz, are available (\cite{wo1}, \cite{wo2}).

A useful point of view comes from  the heuristic formula $A=\c^*(\Gamma )$. Here $\Gamma$ is a discrete quantum group, obtained as a kind of Pontrjagin dual of $G$. Number of discrete group techniques are known to apply to this situation. See e.g. \cite{ve1}, \cite{ve2}.

It is also known from \cite{gr} that quantum permutation groups are in one-to-one correspondence with subalgebras of spin planar algebras constructed in \cite{j1}, \cite{j2}.

Summarizing, a quantum permutation group should be regarded as a mixture of finite, compact and discrete groups, with a flavor of statistical mechanics, knot invariants and planar algebras. Several results are obtained along these lines in \cite{sms}, \cite{gr}, \cite{cl}, \cite{bi1}, \cite{bi2}.

The aim of this work is to bring into the picture some free probability techniques.

The starting point is the classical formula $G(X\ldots X)=G(X)\times_w G(X_n)$ for usual symmetry groups. Here  $X$ is a finite connected graph, $X_n$ is a set having $n$ elements, $X\ldots X$ is the disjoint union of $n$ copies of $X$, and $\times_w$ is a wreath product. A series of free quantum analogues and generalisations, started in \cite{bi2} and continued here, leads to a general formula of type $A(X*Y)=A(X)*_wA(Y)$. Here $X,Y$ are colored graphs, and $*_w$ is a free wreath product.

The corepresentation theory of free wreath products is worked out in two particular situations in \cite{gr}, \cite{bi2}. Our key remark here is that a formula of type $\mu (A*_wB)=\mu(A)\boxtimes\mu(B)$ holds in both cases, where $\mu$ is the associated spectral measure. We conjecture that this formula holds in general, and under mild assumptions on $A$ and $B$.

This is to be related to a planar algebra formula of type $\mu (P*Q)=\mu(P)\boxtimes\mu (Q)$, known to Bisch and Jones (\cite{bj5}). In fact, a general formula of type $A(X*Y)=A(X)*_wA(Y)$, with colored graphs replaced by planar algebras, would be equivalent to the conjecture.

Of particular interest is the case $B=A(X_n)$. Here the conjecture, together with Voiculescu's $S$-transform technique (\cite{v2}) reduces computation of $\mu (X)$ with $X$ homogeneous to that of $\mu(X)$ with $X$ connected and homogeneous. For $n=2$ the conjecture is actually a theorem, and, as an application, we compute $\mu$ for the graph which looks like 2 rectangles. This completes previous classification work for graphs having at most $8$ vertices (\cite{sms}, \cite{gr}).

The paper is organised as follows. 1 is a preliminary section. In 2 and 3 we rearrange some previously known results, with the main improvement that we use spectral measures instead of Poincar\'e series. Among key remarks here is that the passage from classical to quantum permutation groups corresponds to the passage from a Poisson law to a free Poisson law.

In 4 we present a first verification of the conjecture, using free probability tools. In 5 and 6 we prove the formula for graphs, and we deduce from it a second verification of the conjecture, using planar algebra methods. In 7 and 8 we discuss disconnected graphs.

\subsection*{\bf Acknowledgments}
We are grateful to Philippe Biane, Dietmar Bisch and Mireille Capitaine for several useful discussions.

\section{Formalism}

Let $A$ be a unital $\c^*$-algebra. That is, we have an associative algebra with unit $A$ over the field of complex numbers $\c$, with an antilinear antimultiplicative map $a\to a^*$ satisfying $a^{**}=a$, and with a Banach space norm satisfying $||a^*a|| =||a||^2$.

A projection is an element satisfying $p^2=p^*=p$. Two projections are said to be orthogonal when $pq=0$. A partition of the unity is a finite set of projections, which are mutually orthogonal, and sum up to $1$.

\begin{defi}
A magic biunitary matrix is a square matrix $v\in M_n(A)$, all whose rows and columns are partitions of the unity of $A$.
\end{defi}

A magic biunitary is indeed a biunitary, in the sense that both $v$ and its transpose $v^t$ are unitaries. The other word -- magic -- comes from a vague similarity with magic squares.

The basic example comes from the symmetric group $S_n$. Consider the following sets.
$$S_{ij}=\{g\in S_n\mid g (j)=i \}$$

When $i$ is fixed and $j$ varies, or vice versa, these sets form partitions of $S_n$. Thus their characteristic functions $s_{ij}\in\c(S_n)$ form a magic biunitary.

At $n=4$ we have the following key example.
$$d=\begin{pmatrix}p&1-p&0&0\cr 1-p&p&0&0\cr 0&0&q&1-q\cr 0&0&1-q&q \end{pmatrix}$$

Here $p,q$ are projections, say on some Hilbert space $H$. If $p,q$ are chosen to be free, the algebra they generate is isomorphic to $\c^*(\z_2*\z_2)$. This shows that entries of a magic biunitary matrix can generate a non commutative, infinite dimensional algebra.

The universal magic biunitary matrix is constructed by Wang in \cite{wa2}.

\begin{defi}
$A(X_n)$ is the universal $\c^*$-algebra generated by $n^2$ elements $u_{ij}$, with the relations making $(u_{ij})$ a magic biunitary matrix.
\end{defi}

In other words, we have the following universal property. For any magic biunitary matrix $v\in M_n(A)$ there is a morphism of $\c^*$-algebras $A(X_n)\to A$ mapping $u_{ij}\to v_{ij}$.

The matrix $s$ produces a quotient map $A(X_n)\to\c (S_n)$. This is an isomorphism for $n=1,2$, both algebras involved having dimension $1,2$. A similar result holds for $n=3$, because entries of a $3\times 3$ magic biunitary matrix can be shown to commute with each other.

The matrix $d$ produces a quotient map $A(X_4)\to \c^*(\z_2*\z_2)$, which shows that $A(X_4)$ is not commutative, and infinite dimensional. The same holds for $A(X_n)$ with $n\geq 4$.

The algebra $A(X_n)$ has a comultiplication map $\Delta$, making $u$ a corepresentation.
$$\Delta (u_{ij})=\sum u_{ik}\otimes u_{kj}$$

There is also a counit given by $\varepsilon (u_{ij})=\delta_{ij}$, and an antipode given by $S(u_{ij})=u_{ji}$. All three maps are constructed by using the universal property of $A(X_n)$. For instance the counit comes from the fact that the unit matrix $(\delta_{ij})$ is a magic biunitary in $M_n(\c )$.

Thus $A(X_n)$ is a Hopf $\c^*$-algebra in the sense of Woronowicz \cite{wo1}. Associated to it are a compact and a discrete quantum group, according to the following heuristic formula. 
$$A(X_n)=\c (G_n)=\c^*(\Gamma_n)$$

With this formalism, we have an inclusion $S_n\subset G_n$, which for $n=1,2,3$ is an isomorphism. Also, we have a quotient map $\Gamma_4\to\z_2*\z_2$, which shows that $\Gamma_4$ is not abelian, nor finite. Thus its Pontrjagin dual $G_4$ is not a compact group, nor a finite quantum group. The same is true for any $G_n$ with $n\geq 4$, because $G_n$ contains $G_4$ as a quantum subgroup.

The following fundamental result is due to Wang (\cite{wa2}).

\begin{theo}
$A(X_n)$ is the universal Hopf $\c^*$-algebra coacting on $X_n=\{1,\ldots ,n\}$.
\end{theo}

In general, a coaction of a Hopf $\c^*$-algebra $A$ on the set $X_n$ is a morphism of $\c^*$-algebras of the following type, satisfying a coassociativity condition.
$$\alpha :\c (X_n)\to\c (X_n)\otimes A$$

The universal coaction can be defined as a linear map, by the following formula.
$$\alpha (\delta_i )=\sum\delta_j\otimes u_{ji}$$

Since $u$ is a magic biunitary, $\alpha$ is a morphism of $\c^*$-algebras, and since $u$ is a corepresentation, $\alpha$ is coassociative. Universality comes from the fact that a linear map of type $\alpha (\delta_i )=\sum\delta_j\otimes v_{ji}$ is a morphism of $\c^*$-algebras if and only if all rows of $v$ are partitions of unity. By applying the antipode the same must hold for the transpose matrix $v^t$, and this leads to the result.

The coaction $\alpha$ should be interpreted as coming from an action $a(i,\sigma )=\sigma (i)$ of the quantum group $G_n$ on the set $X_n$, via the transposition formula $\alpha\varphi=\varphi a$.
$$a:X_n\times G_n\rightarrow X_n$$

Since $\alpha$ is universal, $a$ is universal as well, so $G_n$ is a kind of quantum analogue of $S_n$.

In what follows, we will be interested in quantum subgroups of $G_n$.

\begin{defi}
A quantum permutation group of $X_n$ corresponds to a pair $(A,v)$, where $A$ is a Hopf $\c^*$-algebra quotient of $A(X_n)$, and $v_{ij}\in A$ is the image of $u_{ij}\in A(X_n)$.
\end{defi}

In other words, we have a Hopf $\c^*$-algebra $A$, and a magic biunitary matrix $v\in M_n(A)$. We assume that $v$ is a corepresentation of $A$, meaning that $\Delta (v_{ij})=\sum v_{ik}\otimes v_{kj}$ and $\varepsilon (v_{ij})=\delta_{ij}$, and that entries of $v$ generate $A$. Observe that the antipode is given by $S(v_{ij})=v_{ji}$.

As with any Hopf $\c^*$-algebra, a main problem regarding $A$ is to find its irreducible corepresentations, and their fusion rules. Since coefficients of $v$ generate $A$, each such corepresentation appears in a tensor power of $v$. So, a first problem is to decompose these tensor powers.
$$v^{\otimes k}=c_k\cdot 1+\sum_{r\neq 1}m_k^r\cdot r$$

The trivial corepresentation $1$ plays here a distinguished role, because computing its multiplicity $c_k$ is the very first question to be asked. By taking characters and by applying the Haar functional we obtain the following formula, where $\chi =v_{11}+\ldots +v_{nn}$ is the character of $v$.
$$\int\chi^k=c_k$$

Consider now a polynomial $\varphi (x)=a_0+a_1x+\ldots +a_kx^k$. By applying it to $\chi$ we get an element $\varphi(\chi)\in A$, then we can apply the Haar functional to this element.
$$\int\varphi (\chi )=a_0+a_1c_1+\ldots +a_kc_k$$

We have here a linear map $\c [X]\to\c$, which extends by density to continuous functions, and is given by integration with respect to a real measure, called spectral measure of $\chi$.

\begin{defi}
The spectral measure of the character of the fundamental corepresentation $\chi =v_{11}+\ldots +v_{nn}$ with respect to the Haar functional is given by
$$\int\varphi (\chi)=\int \varphi (x)\, d\mu (x)$$
for any continuous function $\varphi :\r\to\c$. This is a probability measure on $[0,n]$.
\end{defi}

We will regard $\mu$ as an invariant of $A$, encoding the sequence of numbers $c_k$.

As a first motivating fact, a Kesten type criterion shows that $A$ is amenable in the discrete quantum group sense if and only if $n$ is the upper bound of the support of $\mu$. Also, in certain situations, fusion rules for $A$ can be recovered from $\mu$. See \cite{gr}, \cite{cl}.

Another invariant, equivalent to $\mu$, is the Poincar\'e series of the planar algebra associated to $A$. This planar algebra is the graded union of spaces $P_k$ of fixed points of $v^{\otimes k}$, and its Poincar\'e series, or dimension function, is given by the following formula.
$$f(z)=\sum_{k=0}^\infty c_kz^k=\int\frac{1}{1-z\chi}=\int \frac{1}{1-zx}\, d\mu (x)$$

The Poincar\'e series is the invariant used by Bisch and Jones in their classification program \cite{bj3}, \cite{bj4}. Also, the following related series appears in Jones' fundamental work \cite{j3}.
$$\Theta (q)=q+\frac{1-q}{1+q}\, f\left( \frac{q}{(1+q)^2}\right)$$

This is the generating series of the decomposition of $P$ as sum of planar modules over $TL(n)$. Its coefficients $a_k$ being multiplicities, they satisfy $a_k\geq 0$. These inequalities give a number of conditions on the numbers $c_k$, which can be regarded as fine algebraic properties of $A$.

\section{Poisson laws}

There are a few known computations of spectral measures. These are direct computations of type $A\to\mu (A)$, or decomposition results of type $\mu (A\times B)=\mu (A)\times\mu (B)$.

In this paper we rearrange most of what is known, in the form of two computations, three basic decomposition results, and a conjectural decomposition result.

The two computations are presented in this section. One conclusion will be that, at an asymptotic level, the passage from the classical algebra $\c (S_n)$ to the quantum algebra $A(X_n)$ corresponds to the passage from a Poisson law to a free Poisson law.

Consider first a subgroup $G\subset S_n$. The algebra of complex functions $\c (G)$ is a Hopf $\c^*$-algebra, with comultiplication given by $\Delta (\varphi):(g,h)\to\varphi (gh)$. The corresponding coaction on $X_n$ has as coefficients the characteristic functions of the sets $\{g\in G|g (j)=i\}$.

\begin{prop}
For a subgroup $G\subset S_n$ the spectral measure of $\c (G)$ is
$$\mu =\frac{1}{|G|}\sum_{s=0}^n m_s\,\delta_s$$
where $m_s$ is the number of permutations in $G$ having exactly $s$ fixed points.
\end{prop}

\begin{proof}
We have $\chi (g)=s$, where $s$ is the number of fixed points of $g$. This gives the following formula for the Poincar\'e series, see \cite{gr} for details.
$$f(z)=\frac{1}{|G|}\sum_{s=0}^n\frac{m_s}{1-sz}$$

Since coefficients of $f$ are moments of $\mu$, we conclude that $\mu$ is the spectral measure.
\end{proof}

The Poisson law of parameter $1$ is the following probability measure on the real line.
$$\nu=\frac{1}{e}\sum_{s=0}^\infty \frac{1}{s!}\,\delta_s$$

We can apply Proposition 2.1 to the symmetric group $G=S_n$. The spectral measure appears to be a kind of truncated exponential of $\delta_1-\delta_0$, converging to $\nu$.

\begin{theo}
The spectral measure of $\c (S_n)$ is given by
$$\nu_n =\sum_{t=0}^n\frac{1}{t!}\,(\delta_1-\delta_0)^{*t}$$
where $*$ is the convolution of real measures. In particular we have $\nu_n\to\nu$.
\end{theo}

\begin{proof}
Permutations having exactly $s$ fixed points are obtained by choosing these $s$ points, then by permuting the remaining $l=n-s$ ones in such a way that there is no fixed point. These latter permutations are counted as follows: we start with all permutations, then we substract those having one fixed point, we add those having two fixed points, and so on.
$$m_s=\begin{pmatrix}n\cr s\end{pmatrix}\left( l!-
\begin{pmatrix}l\cr 1\end{pmatrix}(l-1)!+
\begin{pmatrix}l\cr 2\end{pmatrix}(l-2)!-\ldots\right)$$

This allows us to compute the spectral measure.
\begin{eqnarray*}
\nu_n
&=&\frac{1}{n!}\sum_{s=0}^n\begin{pmatrix}n\cr s\end{pmatrix}\left(\sum_{k=0}^{n-s}(-1)^k\frac{(n-s)!}{k!}\right)\,\delta_s\cr
&=&\sum_{s=0}^n\sum_{k=0}^{n-s}(-1)^k\frac{1}{n!}\cdot \frac{n!}{s!(n-s)!}\cdot\frac{(n-s)!}{k!}\,\delta_s\cr
&=&\sum_{s=0}^n\sum_{k=0}^{n-s}\frac{(-1)^k}{s!k!}\,\delta_s\cr
\end{eqnarray*}

We can further improve this formula, by summing over $k$ and over $t=k+s$. 
\begin{eqnarray*}
\nu_n&=&\sum_{t=0}^n\sum_{k=0}^{t}\frac{(-1)^k}{(t-k)!k!}\,\delta_{t-k}\cr
&=&\sum_{t=0}^n\frac{1}{t!}\sum_{k=0}^{t}(-1)^k\begin{pmatrix}t\cr k\end{pmatrix}\,\delta_{t-k}\cr
&=&\sum_{t=0}^n\frac{1}{t!}(\delta_1-\delta_0)^{*t}
\end{eqnarray*}

The last assertion follows by considering the Fourier transform of $\nu_n$.
$$F_n(y)=\sum_{t=0}^n\frac{1}{t!}\left( e^{iy}-1\right)^t$$

With $n\to\infty$ we get the following function.
$$F(y)=e^{( e^{iy}-1)}$$

This is the Fourier transform of the Poisson law, and we are done.
\end{proof}

Another known computation is for the universal algebra $A(X_n)$. Recall first that the normalised semicircular law is the following probability measure on $[-2,2]$.
$$\gamma=\frac{1}{2\pi }\sqrt{4-x^2}\,dx$$

A variable $s$ having this spectral measure is called semicircular. The spectral measure of $s^2$ is called free Poisson law of parameter $1$, and is given by the following formula.
$$\eta =\frac{1}{2\pi}\sqrt{4x^{-1}-1}\,dx$$

The terminology comes from the fact that in the central limit theorem, the Gaussian law from classical probability gets replaced by the semicircular law in free probability. In other words, semicircular means free Gaussian, so square of semicircular means free Poisson.

The measure $\eta$ is also called Marchenko-Pastur law of parameter $1$. See \cite{mp}, \cite{v3}, \cite{vdn}.

\begin{theo}
The spectral measure of $A(X_n)$ is given by $\eta_1=\delta_1$ and
$$\eta_2=\frac{1}{2}(\delta_0+\delta_2)$$
$$\eta_3 =\frac{1}{6}(2\delta_0+3\delta_1+\delta_3)$$
$$\eta_{4+} =\frac{1}{2\pi}\sqrt{4x^{-1}-1}\,dx$$
where $4+$ denotes any positive integer $n\geq 4$. In particular we have $\eta_n\to\eta$.
\end{theo}

\begin{proof}
It is known from \cite{gr} that $A(X_n)$ corresponds to the Temperley-Lieb algebra $TL(n)$, whose Poincar\'e series is computed by counting diagrams, and is given by the following formulae.
$$f_2(z)=\frac{1}{2}\left( 1+\frac{1}{1-2z}\right)$$
$$f_3(z)=\frac{1}{6}\left( 2+\frac{3}{1-z}+\frac{1}{1-3z}\right)$$
$$f_{4+}(z)=\frac{1-\sqrt{1-4z}}{2z}$$

The corresponding measures can be recaptured by using the Stieltjes formula.
$$d\mu (x)=\lim_{t\searrow 0}-\frac{1}{\pi}\,\mbox{Im}\left(G(x+it)\right)\cdot dx$$

Here $G$ is the Cauchy transform of a real measure $\mu$. This is given by $G(\xi )=\xi^{-1}f(\xi^{-1})$, where $f$ is the generating series of the moments of $\mu$. In our situation $f$ is the Poincar\'e series, and we get the following formulae of Cauchy transforms.
$$G_2(\xi)=\frac{1}{2}\left( \frac{1}{\xi}+\frac{1}{\xi-2}\right)$$
$$G_3(\xi)=\frac{1}{6}\left( \frac{2}{\xi}+\frac{3}{\xi -1}+\frac{1}{\xi-3}\right)$$
$$G_{4+}(\xi )=\frac{1-\sqrt{1-4\xi^{-1}}}{2}$$

The Stieltjes formula applies, and gives the measures in the statement.
\end{proof}

This proof, based on counting diagrams, doesn't quite tell us why $\eta_{4+}$ is a free Poisson law, or a Marchenko-Pastur law. In order to get a more enlightening explanation here, some kind of matrix model for $A(X_{4+})$ seems to be needed. So far, the only result in this sense is the one in \cite{cl}, where an explicit realisation of $A(X_4)$ is found. The model there uses a $4\times4$ matrix constructed using quaternions, which shows that $\eta_4$ is indeed the law of the square of a semicircle, appearing as a meridian on the real sphere $S^3$.

\section{Decomposition results}

We discuss now decomposition results of type $\mu (A\times B)=\mu (A)\times\mu (B)$. We have here three basic results, plus a conjecture, unifying two previously known computations from \cite{gr}, \cite{bi2}.

If $A,B$ are given with linear forms $h,k$, the tensor product $A\otimes B$ has a canonical linear form, namely the tensor product $h\otimes k$. Also, the free product $A*B$ has a canonical linear form $h*k$, called free product of $h$ and $k$. See \cite{vdn}.

We regard $A$ and $B$ as being subalgebras of $A\otimes B$ and of $A*B$.

For $x\in A$ and $y\in B$ the spectral measures of $x+y$ and of $xy$ in $A\otimes B$ and in $A*B$ depend only on the spectral measures of $x,y$. The corresponding four operations on real measures are  called additive convolution, multiplicative convolution, free additive convolution, and free multiplicative convolution, and are denoted $*,\times,\boxplus,\boxtimes$.

\begin{defi}
Assume that $x\in A$ and $y\in B$ have spectral measures $\mu$ and $\nu$.

(i) The spectral measure of $x+y\in A\otimes B$ is denoted $\mu *\nu$.

(ii) The spectral measure of $xy\in A\otimes B$ is denoted $\mu\times\nu$.

(iii) The spectral measure of $x+y\in A*B$ is denoted $\mu\boxplus\nu$.

(iv) The spectral measure of $xy\in A*B$ is denoted $\mu\boxtimes\nu$.
\end{defi}

In this definition we assume that $x,y$ are self-adjoint, and so are the variables in (i,ii,iii). As for (iv), the variable $xy$ in there is not exactly self-adjoint, and the spectral measure we refer to is the one of a self-adjoint element having the same moments (for instance $x^{1/2}yx^{1/2}$ when $x>0$, in the tracial case). See Voiculescu's original paper \cite{v2}, or \cite{v3}, \cite{vdn} for details.

Consider now two pairs $(A,u)$ and $(B,v)$ as in Definition 1.3. Denote by $A\times B$ the tensor product, or the free product of $A$ and $B$. This is a Hopf $\c^*$-algebra, and the matrices $u$ and $v$ can be regarded as corepresentations of it. We can form their direct sum and tensor product.
$$u\oplus v=\begin{pmatrix}u&0\cr 0&v\end{pmatrix}$$
$$(u\otimes v)_{ia,jb} =u_{ij}v_{ab}$$

The matrices $u\oplus v$ and $u\otimes v$ are corepresentations of both $A\otimes B$ and $A*B$.  Moreover, $u\oplus v$ is a magic biunitary over $A\otimes B$ and $A*B$, and $u\otimes v$ is a magic biunitary over $A\otimes B$, as one can see by checking the definition. As for $u\otimes v$ over $A*B$, this is not a magic biunitary, for instance because its entries are not self-adjoint.

\begin{prop}
We have the formulae
$$\mu (A\otimes B, u\oplus v)=\mu (A,u)* \mu (B,v)$$
$$\mu (A\otimes B, u\otimes v)=\mu (A,u)\times\mu (B,v)$$
$$\mu (A*B, u\oplus v)=\mu (A,u)\boxplus \mu (B,v)$$
where $*,\times,\boxplus$ are the additive, multiplicative, and free additive convolutions.
\end{prop}

\begin{proof}
The corresponding characters are given by the following formulae.
$$\chi (u\oplus v)=\chi (u)+\chi (v)$$
$$\chi (u\otimes v)=\chi (u)\chi (v)$$

The Haar functional of $A\times B$ being the $\times$ product of those of $A$ and $B$, the variables $\chi(u)\in A\times B$ and $\chi (v)\in A\times B$ are independent or free, and have same spectral measures as $\chi (u)\in A$ and $\chi (v)\in B$. The result follows from the definition of convolutions.
\end{proof}

The free wreath product $A*_wB$ is constructed in \cite{bi1}. This is the algebra generated by $n$ copies of $A$ and a copy of $B$, with the $a$-th copy of $A$ commuting with the $a$-th row of $v$, for any $a$. The maps $\Delta$ and $\varepsilon$ are constructed by using the universal property of $A*_wB$.

\begin{defi}
The free wreath product of $(A,u)$ and $(B,v)$ is given by
$$A*_wB=(A^{*n}*B)/<[u_{ij}^{(a)},v_{ab}]=0>$$
where $n\times n$ is the size of $v$. This has fundamental magic biunitary matrix
$$w_{ia,jb}=u_{ij}^{(a)}v_{ab}$$
and Hopf algebra structure making $w$ a corepresentation.
\end{defi}

There are several approaches to free wreath products, to be discussed in what follows. A very first thing to be noticed is the following diagram.
$$\begin{matrix}
A^{*n}*B&\ &\rightarrow&\ &A*_wB\cr
\ \cr
\downarrow&\ &\ &\ &\downarrow\cr
\ \cr
A*B&\ &\rightarrow&\ &A\otimes B
\end{matrix}$$

Here all maps are surjective morphisms of $\c^*$-algebras, coming from definitions. Vertical maps are obtained by collapsing the $n$ copies of $A$ to a single one, and horizontal maps correspond to certain commutation relations with rows of $v$.

This diagram gives different ways of mixing diagonal elements $u_{ii}$ and $v_{aa}$.
$$\begin{matrix}
\sum u_{ii}^{(a)}v_{aa}&\ &\rightarrow&\ &\sum u_{ii}^{(a)}v_{aa}\cr
\ \cr
\downarrow&\ &\ &\ &\downarrow\cr
\ \cr
\sum u_{ii}v_{aa}&\ &\rightarrow&\ &\sum u_{ii}\otimes v_{aa}
\end{matrix}$$

Here in the upper right corner we have the unknown, namely the character of $w$.

In the lower left corner we have a familiar object, namely the character of $u\otimes v$. This corepresentation $u\otimes v$ is not a magic biunitary, but the corresponding spectral measure can be constructed as in Definition 1.4, and is given by the missing formula in Proposition 3.1:
$$\mu (A*B, u\otimes v)=\mu (A,u)\boxtimes \mu (B,v)$$

The spectral measures of these two elements can be computed in several cases of interest, and turn out to be equal. However, equality doesn't hold in general, and a quite natural assumption here is that all diagonal elements $u_{ii}$, as well as all diagonal elements $v_{aa}$, should have same spectral measure. In terms of Hopf algebras, this leads to the following statement.

\begin{conj}
If $A$ and $B$ have homogeneous skeletons we have the formula
$$\mu (A*_wB,w)=\mu (A,u)\boxtimes \mu(B,v)$$
where $\boxtimes$ is the free multiplicative convolution.
\end{conj}

The notion of skeleton is introduced as follows. The idea is to regard each pair $(B,v)$ as corresponding to the quantum symmetry group of a graph-like object $X=(X_n,D)$. This is definitely possible for all known examples of $(B,v)$, because they appear naturally in this way. In the general case one can take for instance the combinatorial data $D$ to be the planar algebra associated to $(B,v)$, by using the Tannaka-Galois type correspondence in \cite{gr}. 

A skeleton $X=(X_n,D)$ is called homogeneous if for any $i,j\in X_n$ there is a permutation $g\in S_n$ preserving the combinatorial data $D$, such that $g(i)=j$.

All this is probably a bit too heuristic, but finding the precise assumptions in Conjecture 3.1 is rather a matter of proving the conjecture, and this is not among purposes of this paper.

We should mention here that finding a truly delinearised definition for skeletons is quite a subtle problem, beyond the state-of-art of the subject. At level of general planar algebras this is known to be possible, as shown by Kodiyalam and Sunder in \cite{ks}.

In the rest of this section we will just explain why some assumptions are needed, and why these should correspond to homogeneity of skeletons of $A$ and $B$.

Our first task is to construct a counterexample. We use the following simple fact.

\begin{prop}
We have a canonical isomorphism
$$(A,u)*_w(B_1*B_2,v_1\oplus v_2)=(A*_wB_1) *(A*_wB_2)$$
for any three pairs $(A,u)$, $(B_1,v_1)$ and $(B_2,v_2 )$.
\end{prop}

\begin{proof}
The algebra on the left is generated by $n_1+n_2$ copies of $A$ and a copy of $B_1*B_2$, subject to certain commutation relations. These relations say that the first $n_1$ copies of $A$ commute with corresponding rows of the matrix $v_1$, and that the remaining $n_2$ copies of $A$ commute with corresponding rows of the matrix $v_2$. In other words, when moving the copy of $B_2$ from the $n_1+n_2+2$-th position to the $n_1+1$-th position, we get the algebra on the right.
\end{proof}

Now if we assume that Conjecture 3.1 holds without assumptions on $(B,v)$, we would get from Propositions 3.1 and 3.2 a kind of distributivity formula for $\boxtimes$ with respect to $\boxplus$.
$$\mu\boxtimes (\mu_1\boxplus\mu_2)=(\mu\boxtimes\mu_1)\boxplus (\mu\boxtimes\mu_2)$$

This formula is known to be wrong, so some assumptions on $(B,v)$ are needed in Conjecture 3.1. We believe that homogeneity of the skeleton, not satisfied by $(B_1*B_2,v_1\oplus v_2)$, but satisfied in the two cases where we know how to prove the conjecture, is the good one.

As for the assumption on $(A,u)$, we are less confident about it. Let us mention here that Conjecture 3.1 should be regarded as an analogue of the planar algebra formula $\mu (P*Q)=\mu (P)\boxtimes\mu (Q)$ of Bisch and Jones (\cite{bj5}). In fact, we have the following conjectural computation.
\begin{eqnarray*}
\mu (A*_wB)
&=&\mu (P(A*_wB))\cr
&=&\mu (P(A)*P(B))\cr
&=&\mu (P(A))\boxtimes \mu(P(B))\cr
&=&\mu (A)\boxtimes \mu(B)
\end{eqnarray*}

Here the first and fourth equality come from the Tannaka-Galois type correspondence $A\to P(A)$ from \cite{gr}. The third equality comes from the formula of Bisch and Jones. The second equality comes from the following conjectural formula.
$$P(A*_wB)=P(A)*P(B)$$

This formula should be regarded as a version of Conjecture 3.1. The above computation shows that the formula is stronger than the conjecture, but one can probably go in the other sense as well, by using the following basic lemma.

\begin{lemm}
Assume that $(A,u)$ and $(B,v)$ correspond to quantum permutation groups of $X_n$, and let $f:A\to B$ be a morphism of $\c^*$-algebras such that $f(u_{ij})=v_{ij}$.

(i) $f$ is a surjective morphism of Hopf $\c^*$-algebras.

(ii) The $k$-th moment of $\mu (A)$ is smaller than the $k$-th moment of $\mu (B)$, for any $k$.

(iii) $f$ is an isomorphism if and only if $\mu (A)=\mu (B)$. 
\end{lemm}

\begin{proof}
The first two assertions follows from definitions. The third one is well-known, and follows by taking a basis of coefficients of irreducible corepresentations, as in \cite{wo2}. 
\end{proof}

As a conclusion, what is behind the conjecture is that the free wreath product operation $*_w$ should correspond to a kind of free product operation $*$ at level of skeletons. This approach suggests that $(A,u)$ and $(B,v)$ will have to play at some point symmetric roles, and this is why homogeneity of the skeleton of $(A,u)$ is included in Conjecture 3.1.

\section{Finite groups}

The simplest example of a free wreath product is with $B=\c (G)$, where $G$ is a finite group acting on itself. That is, we consider $G$ as being a subgroup of $S_n$, where $n$ is the order of $G$. The matrix $v$ has a particularly simple form.
$$v=\begin{pmatrix}
\delta_1&\delta_{12}&\ldots &\delta_{1n}\cr
\delta_{21}&\delta_1&\ldots &\delta_{2n}\cr
\dots&\dots&\dots&\dots\cr
\delta_{n1}&\delta_{n2}&\ldots&\delta_1
\end{pmatrix}$$

Here $\delta_1$ is the Dirac mass at the unit element $1\in G$, and for $i,j\in X_n$ we denote by $\delta_{ij}$ the Dirac mass at the unique element $g\in G$ satisfying $g(j)=i$. We see that each row of $v$ consists of Dirac masses at elements of $G$, so in particular its entries generate $\c (G)$. This shows that a free wreath product by $\c (G)$ has the following decomposition.
$$A*_w\c (G)=A^{*n}\otimes\c (G)$$

This is pointed out in \cite{bi2}, where the case of $G=\z_2$ is worked out explicitely, with a complete discussion of corepresentations of the free wreath product. The fusion rules found there allow one to compute the spectral measure of the free wreath product. In what follows we present this computation of spectral measure, along with a generalisation to arbitrary groups $G$.

We use Voiculescu's $S$-transform, whose log linearises $\boxtimes$.
$$S(\mu\boxtimes\nu)=S(\mu)S(\nu)$$

The $S$-transform of a measure $\mu$ is constructed by introducing series $\psi,\chi$ as follows (\cite{v2}).
$$\psi (z)=f(z)-1\hskip 1cm \chi\circ\psi =\psi\circ\chi=id\hskip 1cm S(z)=\frac{1+z}{z}\,\chi (z)$$

Here $f$ is the series having as coefficients the moments of $\mu$. In case $\mu$ is the spectral measure associated to a Hopf $\c^*$-algebra $A$, this series $f$ is the Poincar\'e series of $A$.

It is useful to compute both the spectral measure of $\c (G)$, and its $S$-transform.

\begin{prop}
For a group $G$ acting on itself the corresponding spectral measure of $\c (G)$ and its $S$-transform are given by
$$\mu=\frac{n-1}{n}\,\delta_0+\frac{1}{n}\,\delta_n$$
$$S(z)=\frac{1+z}{1+nz}$$
where $n$ is the number of elements of $G$.
\end{prop}

\begin{proof}
The identity of $G$ has $n$ fixed points, and the other $n-1$ elements, none. Together with Proposition 2.1 this gives the formula of $\mu$, then we get $f,\psi ,\chi$.
\begin{eqnarray*}
\mu=\frac{n-1}{n}\,\delta_0+\frac{1}{n}\,\delta_n
&\Rightarrow&f(z)=1+\frac{z}{1-nz}\cr
&\Rightarrow&\psi (z)=\frac{z}{1-nz}\cr
&\Rightarrow&z=\frac{\chi(z)}{1-n\chi(z)}\cr
&\Rightarrow&\chi (z)=\frac{z}{1+nz}
\end{eqnarray*}

Multiplying by $1+z^{-1}$ gives the formula in the statement.
\end{proof}

\begin{theo}
If $G$ is a finite group acting on itself then
$$\mu (A*_w\c(G))=\mu (A)\boxtimes \mu(\c (G))$$
where $\boxtimes$ is the free multiplicative convolution.
\end{theo}

\begin{proof}
We denote by $\chi, \chi_l$ the characters of $u,w$, and by $\mu,\mu_l$ the associated spectral measures. The tensor product decomposition of the free wreath product gives the following formula.
$$\chi_l =\left(\chi^{(1)}+\ldots +\chi^{(n)}\right)\otimes\delta_1$$

We raise both terms to the power $k\geq 1$. 
$$\chi_l^k =\left(\chi^{(1)}+\ldots +\chi^{(n)}\right)^k\otimes\delta_1$$

The Haar functional of the tensor product being the tensor product of Haar functionals, we get the following formula.
$$\int \chi_l^k=\frac{1}{n}\,\int\left(\chi^{(1)}+\ldots +\chi^{(n)}\right)^k$$

We get a formula which is valid for any polynomial $\varphi$.
$$\int \varphi (\chi_l )=\frac{n-1}{n}\,\varphi (0)+\frac{1}{n}\,\int\varphi\left(\chi^{(1)}+\ldots +\chi^{(n)}\right)$$

Now remember that elements $\chi^{(a)}$ are free, and each of them has spectral measure $\mu$. This gives the following identity.
$$\int\varphi (x)\,d\mu_l(x)=\frac{n-1}{n}\,\varphi (0)+\frac{1}{n}\,\int \varphi (x)\,d\mu^{\boxplus n}(x)$$

This finishes computation of the spectral measure on the left.
$$\mu_l=\frac{n-1}{n}\,\delta_0+\frac{1}{n}\,\mu^{\boxplus n}$$

The spectral measure on the right is given by the following formula.
$$\mu_r=\mu\boxtimes \left( \frac{n-1}{n}\,\delta_0+\frac{1}{n}\,\delta_n\right)$$

We can use here the following general formula, discussed below.
$$\frac{n-1}{n}\,\delta_0+\frac{1}{n}\,\mu^{\boxplus n}=\mu\boxtimes \left( \frac{n-1}{n}\,\delta_0+\frac{1}{n}\,\delta_n\right)$$

Together with the above descriptions of $\mu_l$ and $\mu_r$, this finishes the proof.
\end{proof}

The free probability formula used in the end of proof of Theorem 4.1 is proved by Nica and Speicher in \cite{ns} by using non-crossing partitions. A second proof is given by Voiculescu in the appendix of \cite{ns}, by using random matrices.

The proof below was already written when we learned about \cite{ns}, and this is the main reason why we include it. The other reason is that we hope that a suitable analytic function procedure can be used in order to simplify the operation $A\to \mu (A)$, so a complete proof of Theorem 4.1 using analytic functions is probably useful.

We use Voiculescu's $R$-transform, which linearises $\boxplus$. 
$$R(\mu\boxplus\nu)=R(\mu)+R(\nu)$$

The $R$-transform of a measure $\mu$ is defined by introducing series $G,K$ as follows (\cite{v1}).
$$G(\xi)=\xi^{-1}f(\xi^{-1})\hskip 1cm K\circ G=G\circ K=id\hskip 1cm R(z)=K(z)-\frac{1}{z}$$

Here, as usual, $f$ denotes the series having as coefficients the moments of $\mu$.

A first remark is that both $R$ and $S$ are given, up to some normalisations, by an inversion of series. We have the following formulae connecting $R$ and $S$.

\begin{lemm}
We have the formulae
$$R(zS(z))=\frac{1}{S(z)}\hskip 2cm S(zR(z))=\frac{1}{R(z)}$$
where $R$ and $S$ are the $R$ and $S$ transforms of the same measure $\mu$.
\end{lemm}

\begin{proof}
The starting point is the following formula, relating $G$ to $\psi$.
$$G(\xi)=\frac{1}{\xi}\,f\left(\frac{1}{\xi}\right)=\frac{1}{\xi}+\frac{1}{\xi}\,\psi\left(\frac{1}{\xi}\right)$$

We make the replacement $\xi\to 1/\chi(z)$.
$$G\left(\frac{1}{\chi (z)}\right)=\chi(z)+\chi(z)z$$

By applying $K$ we get a formula relating $K$ and $\chi$.
$$\frac{1}{\chi (z)}=K(\chi(z)(1+z))$$

In terms of $R$ and $S$, we get the following equality.
$$\frac{1+z}{z}\cdot \frac{1}{S(z)}=K(zS(z))=R(zS(z))+\frac{1}{zS(z)}$$

This gives the first formula. For the second one, we start again with the equality relating $G$ to $\psi$, and we make the replacement $\xi\to K(z)$.
$$z=\frac{1}{K(z)}+\frac{1}{K(z)}\psi\left(\frac{1}{K(z)}\right)$$

This can be written in the following way.
$$\psi\left(\frac{1}{K(z)}\right) =zK(z)-1$$

By applying $\chi$ we get a second formula relating $K$ and $\chi$.
$$\frac{1}{K(z)}=\chi(zK(z)-1)$$

On the other hand, we have the following equality.
$$S(zR(z))=\frac{1+zR(z)}{zR(z)}\cdot \chi (zR(z))=\frac{K(z)}{R(z)}\cdot \chi (zR(z))$$

Together with $zR(z)=zK(z)-1$, this gives the second formula.
\end{proof}

\begin{theo}
We have the Nica-Speicher-Voiculescu formula
$$\frac{n-1}{n}\,\delta_0+\frac{1}{n}\,\mu^{\boxplus n}=\mu\boxtimes \left( \frac{n-1}{n}\,\delta_0+\frac{1}{n}\,\delta_n\right)$$
for any probability measure $\mu$ on the real line.
\end{theo}

\begin{proof}
We denote by $\mu_l$ and $\mu_r$ the measures on the left and on the right. Let also $\mu_+$ be the measure $\mu^{\boxplus n}$, and $\mu_\delta$ be the measure in the bracket.

For any measure of type $\mu_x$ we denote by $f_x,G_x,K_x,R_x,\psi_x,\chi_x,S_x$ the series involved in the construction of $R$ and $S$ transforms.

The starting point is the following equality, coming from the fact that $R$ linearises $\boxplus$.
$$R_+(z)=nR(z)$$

Together with the above lemma, this gives the following formula.
\begin{eqnarray*}
S_+(nz)
&=&S_+(nzS(z)R(zS(z)))\cr
&=&S_+(zS(z)R_+(zS(z)))\cr
&=&\frac{1}{R_+(zS(z))}\cr
&=&\frac{1}{nR(zS(z))}\cr
&=&\frac{S(z)}{n}
\end{eqnarray*}

We use now the following identity.
$$\psi_+(z)=f_+(z)-1=n(f_l(z)-1)=n\psi_l(z)$$

At level of $\chi$ functions, this translates as follows.
$$\chi_l(z)=\chi_+\psi_+\chi_l(z)=\chi_+(n\psi_l\chi_l(z))=\chi_+(nz)$$

We have now all ingredients for computing $S_l$.
\begin{eqnarray*}
S_l(z)
&=&\frac{1+z}{z}\,\chi_l(z)\cr
&=&\frac{1+z}{z}\,\chi_+(nz)\cr
&=&\frac{1+z}{z}\cdot\frac{nz}{1+nz}\, S_+(nz)\cr
&=&\frac{(1+z)n}{1+nz}\cdot\frac{S(z)}{n}\cr
&=&\frac{1+z}{1+nz}\, S(z)
\end{eqnarray*}

Proposition 4.1 shows that the right term is $S_\delta(z)S(z)=S_r(z)$, and we are done.
\end{proof}

\section{Colored graphs}

Just as in the classical case, a natural way to get quantum permutation groups is to consider simple algebraic structures like graphs and construct their quantum automorphism groups.

Given a finite graph $X$, the idea is to label its vertices $1,\ldots ,n$, then to consider an algebra of form $A(X)=A(X_n)/J$, where $J$ is an ideal expressing preservation of edges.

This will be explained later on. For the moment we have to fix some notations regarding $X$. In previous papers \cite{sms}, \cite{gr}, \cite{bi1}, \cite{bi2} we use several kinds of finite graphs - with edges oriented or not, colored or not - plus finite metric spaces.

For the purposes of this paper, most convenient is to start with the following definition.

\begin{defi}
In this paper a colored graph $X=(V,E)$ is a finite set $V=X_n$ together with a partition $E$ of the form
$$(V\times V)-\Delta_V=E_1\sqcup\ldots\sqcup E_p$$
where $\Delta_V$ is the diagonal of $V\times V$.
\end{defi}

The terminology is of course not standard, we just use it in order to simplify presentation. In fact colors are obviously missing, so $X$ should be rather called ``colorable'' graph. We use ``colored'' because most interesting examples come along with a natural coloring.

In general, we can use a color set $C$, and an injective map $c$ assigning colors to edges.
$$c:\{1,\ldots ,p\}\to C$$

We get a picture of $X$ in the following way. We call points elements of $V$, we draw them in the plane, then we draw an oriented edge $\vec{ij}$ between any two points $i\neq j$ and we color it $c(k)$, with $k$ given by $(i,j)\in E_k$. This kind of picture is to be called colored graph as well.

As a first example, consider a finite metric space $X$. This can be regarded as a colored graph, with edges colored by their lenghts. For instance when all distances are equal, say to $a\in (0,\infty)$, we get a monocolored graph, called simplex and denoted $X_n^a$.

It is convenient to call simplex any monocolored graph. In case the unique color is already specified, say it is the color $a$, the simplex is denoted $X_n^a$. In case the unique color is not specified, our favorite choice will be the real number $a=0$.

The second example is with an oriented graph $X=(V,E)$. Here $V=X_n$ is a finite set, called set of vertices, and $E$ is a set of oriented edges between pairs of distinct vertices.
$$E\subset V\times V-\Delta_V$$

We get a partition in the following way, where $\bar{E}$ is the set of missing edges.
$$V\times V-\Delta_V=E\sqcup\bar{E}$$

We can regard $X$ as a colored graph, with color $1$ assigned to edges, and color $0$ assigned to missing edges. For instance $E=\emptyset$ gives the simplex $X_n^0$ and $\bar{E}=\emptyset$ gives the simplex $X_n^1$.

In both the metric space and the graph example colors are complex numbers, and the matrix $d_{ij}=c(k)$ is the Laplacian of $X$. This suggests the following definition.

\begin{defi}
The Laplacian of $X$ is the matrix given by $d_{ii}=0$ and
$$d_{ij}=c(k)$$
with $k$ given by $(i,j)\in E_k$, where $c:\{1,\ldots ,p\}\to\c$ is a fixed injective map.
\end{defi}

If $X$ comes from a metric space or from an oriented graph, or more generally if $X$ comes along with a natural complex coloring map $c$, this map will be the one used for constructing the Laplacian, unless otherwise specified.

The fact that the Laplacian depends on the choice of the complex coloring map $c$ is not a problem, because all constructions in this paper do not depend on this choice. Actually, one of the reasons of presenting things this way is that we want to use independence from the choice of colors, for instance by using our favorite color $a=0$, whenever needed.

We will be interested in commutation relations of type $vd=dv$, with $v$ magic biunitary matrix. The Stone-Weierstrass theorem shows that validity of such a commutation relation doesn't depend on the choice of $c$. See Theorem 2.2 in \cite{gr}.

In other words, we can consider the ideal $J_{min}\subset A(X_n)$ generated by the relations $vd=dv$. Consider also the ideal $J_{com}$ generated by all commutation relations between $v_{ij}$'s.

\begin{prop}
We have a canonical isomorphism
$$A(X_n)/J_{max}\simeq \c (G)$$
where $J_{max}=<J_{min},J_{com}>$, and where $G$ is the symmetry group of $X$.
\end{prop}

\begin{proof}
The quotient $A(X_n)/J_{com}$ is canonically isomorphic to $\c (S_n)$, then further quotienting by relations $[v,d]=0$ gives $\c (G)$. See \cite{gr} for details.
\end{proof}

This result suggests to define $A(X)$ to be the quotient of $A(X_n)$ by an ideal $J$ satisfying $J_{min}\subset J\subset J_{max}$. The choice $J=J_{max}$ is of course not very interesting. Another natural choice is $J=J_{min}$, used in \cite{sms}, \cite{gr}. An intermediate choice is made in \cite{bi1}, \cite{bi2}.

In this paper we extend and generalise results in \cite{bi2} to algebras defined using $J_{min}$, in order to get a general formula of type $A(X*Y)=A(X)*_wA(Y)$. Together with results in \cite{gr}, this will lead to a second verification of Conjecture 3.1. 

\begin{defi}
The algebra $A(X)$ is the quotient of $A(X_n)$ by the relations $vd=dv$.
\end{defi}

The pair $(A(X),v)$ satisfies conditions in Definition 1.3. In view of Proposition 5.1, the underlying quantum permutation group is an analogue of the automorphism group of $X$.

As an example, for a simplex we can take the Laplacian to be $d=0$, so we get $A(X_n)$ itself. We can see here the difference with the formalism in \cite{bi1}, \cite{bi2}, where the graph having $n$ vertices and no edges gives a different $A$ algebra from that of the graph having $n$ vertices and all possible edges. In fact, all definitions in this section are motivated by the present choice of $A(X)$.

An interesting situation is when $X$ is the $n$-gon. For any $n\neq 4$ we have $A(X)=\c (D_n)$, but for $n=4$ the algebra $A(X)$ is not commutative, and infinite dimensional. See \cite{gr}.

The relations defining $A(X)$ might be rewritten in several convenient manners.

\begin{prop}
The algebra $A(X)$ is the quotient of $A(X_n)$ by the relations
$$\sum_{k,(k,j)\in E_r}v_{ik}=\sum_{k,(i,k)\in E_r}v_{kj}$$
with the number $r$ ranging over the set $\{ 1,\ldots ,p\}$.
\end{prop}

\begin{proof}
We choose a Laplacian $d$, corresponding to a complex coloring map $c$. If $d_r$ denotes the Laplacian of the oriented graph $X_r=(V,E_r)$, we have the following formula.
$$d=\sum_{r=1}^p c(r)d_r$$

The matrix $d_r$ multiplies by $v$ in the following way.
$$(vd_r)_{ij}=\sum_kv_{ik}(d_r)_{kj}=\sum_{k,(k,j)\in E_r}v_{ik}$$
$$(d_rv)_{ij}=\sum_k(d_r)_{ik}v_{kj}=\sum_{k,(i,k)\in E_r}v_{kj}$$

This shows that the $r$-th relation in the statement is equivalent to the commutation relation $vd_r=d_rv$. On the other hand independence of $A(X)$ from the choice of the Laplacian shows that $vd=dv$ is equivalent to $vd_r=d_rv$ for any $r$, and this gives the result. 
\end{proof}

Given two colored graphs $X,Y$, we can consider the colored graph $X*Y$ obtained by putting a copy of $X$ at each vertex of $Y$. In other words, we have the following definition.

\begin{defi}
Let $X=(T,E)$ and $Y=(Z,F)$ be two colored graphs, 
with $E=\{E_1,\ldots , E_p\}$ and 
$F= \{F_1 \ldots , F_q\}$. 
We define subsets of $T\times Z$ times itself in the following way.
$$E_r^\circ=\{(i\alpha,j\alpha)\mid (i,j)\in E_r,\,\alpha\in Z\}$$ 
$$F_s^\circ=\{(i\alpha,j\beta)\mid i,j\in T,\,(\alpha,\beta)\in F_s\}$$

These determine a colored graph $X*Y$, called free product of $X$ and $Y$.
\end{defi}

The terminology is of course not standard, the free product of colored graphs not being  a coproduct in a categorical sense. We call it like this because it is expected to be compatible with the planar algebra free product $*$ of Bisch and Jones \cite{bj5}.

As a first example, consider a free product of $s$ simplices.
$$Z=Y_1^{\varepsilon_1}*\ldots *Y_s^{\varepsilon_s}$$

As pointed out in \cite{gr}, we have a nice interpretation of $Z$ if we choose the corresponding $s$ colors to be an increasing sequence of infinitesimals.
$$\varepsilon_1<<\varepsilon_2<<\ldots <<\varepsilon_s$$

The free product is obtained by starting with $Y_s$, then putting a copy of $Y_{s-1}$ at each vertex of $Y_s$, a copy of $Y_{s-2}$ at each vertex of $Y_{s-1}$, and so on until $Y_1$. What we get is a kind of metric space, by using the infinitesimal summing conventions $\varepsilon_i+\varepsilon_j=\varepsilon_i$ for $i>j$.

A simplex is called generic if it has $n\geq 4$ vertices. The terminology comes from the proof below, where numbers $n$ become Jones indices, called generic when bigger than $4$.

\begin{theo}
If $Z$ is a free product of $s$ generic simplices then
$$\mu (A(Z))=\eta^{\boxtimes s}$$
where $\eta$ is the free Poisson law of parameter $1$.
\end{theo}

\begin{proof}
It is shown in \cite{gr} that the corresponding $s$ Laplacians satisfy Landau's exchange relations in \cite{la}. This shows that the planar algebra associated to $A(Z)$ is a Fuss-Catalan algebra on $s$ colors, whose Poincar\'e series is computed in the generic case by Bisch and Jones in \cite{bj1}.
$$f(z)=\sum_{k=0}^{\infty}\frac{1}{sk+1}\begin{pmatrix}(s+1)k\cr k\end{pmatrix}z^k$$

As pointed out in \cite{bj2}, this series is a solution of the following  equation.
$$f(z)=1+zf(z)^{s+1}$$

This can be used for computing the $S$-transform of the corresponding spectral measure.
\begin{eqnarray*}
f(z)=1+zf(z)^{s+1}
&\Rightarrow&\psi (z)=z(1+\psi (z))^{s+1}\cr 
&\Rightarrow&z=\chi (z)(1+z)^{s+1}\cr
&\Rightarrow&\chi (z)=\frac{z}{(1+z)^{s+1}}\cr
&\Rightarrow&S(z)=\frac{1}{(1+z)^{s}}
\end{eqnarray*}

With $s=1$ this gives the $S$-transform of $\eta_{4+}=\eta$. Now back to arbitrary $s$, we get the following equality.
$$S (z)=\left(\frac{1}{1+z}\right)^s=(S\eta (z))^s$$

The result follows from the fact that $\log S$ linearises $\boxtimes$. 
\end{proof}

\section{Free products}

The purpose of this section is to establish the formula $A(X*Y)=A(X)*_wA(Y)$. This will be done via a modification plus generalisation of a proof in \cite{bi2}, where free wreath products are constructed, and where a first such decomposition result is found.

The idea is to construct a pair of inverse morphisms, by using universal properties of various algebras involved. The morphism from left to right is easy to construct, and this is done in proof of Theorem 6.1 below. In the other sense, we have first to construct inside $A(X*Y)$ analogues $U,V$ of matrices $u,v$, and this is the purpose of next two lemmas.

The magic biunitary matrix associated to $A(X*Y)$ is denoted $W_{i\alpha,j\beta}$. Here double indices $i\alpha ,j\beta\in T\times Z$ are produced by indices $i,j\in T$ and $\alpha ,\beta\in Z$, coming from Definition 5.4.

\begin{lemm}
We can define an element $V_{\alpha\beta}$ by the formula
$$V_{\alpha\beta}=\sum_kW_{i\alpha,k\beta}=\sum_kW_{k\alpha,j\beta}$$
for any $i,j$, and the resulting matrix $V$ is a magic biunitary.
\end{lemm}

\begin{proof}
The relations in Proposition 5.2 are written in the following way.
$$\sum_{k,(k,j)\in E_r}W_{i\alpha,k\beta}=\sum_{k,(i,k)\in E_r}W_{k\alpha,j\beta}$$
$$\sum_{k,\gamma,(\gamma,\beta)\in F_s}W_{i\alpha,k\gamma}=\sum_{k,\gamma,(\alpha,\gamma)\in F_s}W_{k\gamma,j\beta}$$

The equality in the statement is obtained from these relations.
$$\sum_k W_{i\alpha,k\beta} =
\sum_r \sum_{k,(k,j) \in E_r} W_{i\alpha,k\beta} =
\sum_r \sum_{k, (i,k) \in E_r} W_{k\alpha,j\beta} =
\sum_k W_{k\alpha,j\beta}$$

The entries of $V$ are easily seen to be self-adjoint.
$$V_{\alpha\beta}^*=\sum_i\left(W_{i\alpha,j\beta}\right)^*=\sum_i W_{i\alpha,j\beta}=V_{\alpha\beta}$$

We check now the conditions regarding sums on rows and columns.
$$\sum_\alpha V_{\alpha\beta}=\sum_{i\alpha} W_{i\alpha,j\beta}=1$$
$$\sum_\beta V_{\alpha\beta}=\sum_{j\beta} W_{i\alpha,j\beta}=1$$

The fact that rows of $V$ are partitions of unity is proved as follows.
\begin{eqnarray*}
V_{\alpha\beta}V_{\alpha\gamma}
&=&\left( \sum_jW_{i\alpha,j\beta}\right)\left(\sum_kW_{i\alpha,k\gamma}\right)\cr
&=&\sum_{jk} W_{i\alpha,j\beta}W_{i\alpha,k\gamma}\cr
&=&\delta_{\beta\gamma}\sum_jW_{i\alpha,j\beta}\cr
&=&\delta_{\beta\gamma}V_{\alpha\beta}
\end{eqnarray*}

A similar computation gives the following formula, valid for any $\alpha\neq\gamma$.
$$V_{\alpha\beta}V_{\gamma\beta}=\left( \sum_iW_{i\alpha,j\beta}\right)\left(\sum_kW_{k\gamma,j\beta}\right) =\sum_{ik} W_{i\alpha,j\beta}W_{k\gamma,j\beta}=0$$

Thus entries of $V$ are orthogonal on each column, and this finishes the proof.
\end{proof}

\begin{lemm}
For any $\alpha$ the matrix
$$U_{ij}^\alpha =\sum_\beta W_{i\alpha,j\beta}$$
is a magic biunitary.
\end{lemm}

\begin{proof}
The entries of $U$ are easily seen to be self-adjoint.
$$\left(U_{ij}^\alpha\right)^*=\sum_\beta\left(W_{i\alpha,j\beta}\right)^*=\sum_\beta W_{i\alpha,j\beta}=U_{ij}^\alpha$$

We check now the conditions regarding sums on rows and columns.
$$\sum_jU_{ij}^\alpha=\sum_{j\beta}W_{i\alpha,j\beta}=1$$
$$\sum_iU_{ij}^\alpha=\sum_{i\beta}W_{i\alpha,j\beta}=\sum_\beta V_{\alpha\beta}=1$$

From the fact that $W$ is magic we get that all rows of $U$ are partitions of unity.
$$U_{ij}^\alpha U_{ik}^\alpha = \sum_{\beta\gamma}
W_{i\alpha,j\beta} W_{i\alpha,k\gamma}= 
\delta_{jk}  \sum_{\beta}
W_{i\alpha,j\beta} =\delta_{jk}U_{ij}^\alpha$$

Also, since $V$ is magic we get the following identity, valid for any $\beta\neq\gamma$.
\begin{eqnarray*}
W_{i\alpha,j\beta} W_{k\alpha,j\gamma}
&=&W_{i\alpha,j\beta}\left( \sum_rW_{r\alpha,j\beta}\right)\left( \sum_sW_{s\alpha,j\gamma}\right) W_{k\alpha,j\gamma}\cr
&=&W_{i\alpha,j\beta}V_{\alpha\beta}V_{\alpha\gamma}W_{k\alpha,j\gamma}\cr
&=&0
\end{eqnarray*}

This gives the following formula, valid for any $i\neq k$.
$$U_{ij}^\alpha U_{kj}^\alpha =\sum_{\beta\gamma}
W_{i\alpha,j\beta} W_{k\alpha,j\gamma}=\sum_\beta W_{i\alpha,j\beta} W_{k\alpha,j\beta}=0$$

Thus entries of $U$ are orthogonal on each column, and this finishes the proof.
\end{proof}

\begin{theo}
We have a canonical isomorphism
$$A(X*Y)=A(X)*_wA(Y)$$
for any two colored graphs $X$ and $Y$.
\end{theo}

\begin{proof}
We use the first assertion in Lemma 3.1. What is left to do is to construct morphisms of $\c^*$-algebras from left to right and from right to left, making $W_{i\alpha ,j\beta}$ correspond to $w_{i\alpha,j\beta}$. These two morphisms will be inverse Hopf $\c^*$-algebra isomorphisms.

``$\rightarrow$''. The matrix $w$ being a magic biunitary, we just have to prove that it commutes with the Laplacian of $X*Y$. For, we use relations in Proposition 5.4, which for $X*Y$ correspond to the first two formulae in proof of Lemma 6.1. The first one is checked as follows.
$$\sum_{k,(k,j)\in E_r}w_{i\alpha,k\beta}=\sum_{k,(k,j)\in E_r}u_{ik}^{(\alpha)}v_{\alpha\beta}=\sum_{k,(i,k)\in E_r}u_{kj}^{(\alpha)}v_{\alpha\beta}=\sum_{k,(i,k) \in E_r} w_{k\alpha,j\beta}$$

As for the second relation, this is checked as follows.
$$\sum_{k,\gamma,(\gamma,\beta)\in F_s}w_{i\alpha,k\gamma}=\sum_{k,\gamma,(\gamma,\beta)\in F_s}u_{ik}^{(\alpha)}v_{\alpha\gamma}=\sum_{\gamma,(\gamma,\beta)\in F_s}v_{\alpha\gamma}$$
$$=\sum_{\gamma,(\alpha,\gamma)\in F_s}v_{\gamma\beta} =
\sum_{k,\gamma,(\alpha,\gamma)\in F_s} u_{kj}^{(\gamma)}v_{\gamma\beta}=\sum_{k,\gamma,(\alpha,\gamma)\in F_s}w_{k\gamma,j\beta}$$

``$\leftarrow$'' The matrices $U,V$ being magic biunitaries, we just have to prove that they satisfy the same relations as $u,v$. First is commutation of $U$ with the Laplacian of $X$.
$$\sum_{k, (k,j)\in E_r}U_{ik}^\alpha =
\sum_{\beta,k,(k,j)\in E_r} W_{i\alpha,k\beta} = 
\sum_{\beta, k, (i,k) \in E_r} W_{k\alpha,j\beta}
=\sum_{k, (i,k) \in E_r}U_{kj}^\alpha$$

Second comes commutation of $V$ with the Laplacian of $Y$.
\begin{eqnarray*}
\sum_{\gamma, (\gamma,\beta)\in F_s}V_{\alpha\gamma}
&=&
\sum_k \sum_{\gamma, (\gamma,\beta) \in F_s} W_{i\alpha,k\gamma}\cr
&=&\sum_k \sum_{\gamma, (\alpha,\gamma) \in F_s} W_{k\gamma,j\beta}\cr
&=&\sum_{\gamma, (\alpha, \gamma) \in F_s} V_{\gamma\beta}
\end{eqnarray*}

Finally, the commutation condition defining free wreath products is checked as follows.
$$U_{ij}^{\alpha}V_{\alpha\beta} =
\sum_{\gamma} W_{i\alpha,j\gamma} V_{\alpha\beta}= 
\sum_{\gamma k} W_{i\alpha,j\gamma} W_{i \alpha,k \beta}
= W_{i \alpha,j\beta}$$
$$V_{\alpha\beta}U_{ij}^{\alpha}= V_{\alpha\beta} \sum_{\gamma}W_{i\alpha,j\gamma}=
\sum_{\gamma k} W_{i\alpha,k\beta} W_{i \alpha,j\gamma}=W_{i \alpha,j\beta}$$

Thus we have an arrow from right to left, and this finishes the proof.
\end{proof}

As a first application, we get a second verification of Conjecture 3.1.

\begin{coro}
We have the formula
$$\mu (A(X)*_wA(Y))=\mu (A(X))\boxtimes\mu (A(Y))$$
for any two free products of generic simplices $X$ and $Y$.
\end{coro}

\begin{proof}
This follows by combining Theorems 5.1 and 6.1.
\end{proof}

The genericity assumption in this statement can actually be removed. One way is by using the above proof, with the formulae for generic indices in Theorem 5.1 and its proof replaced by formulae for arbitrary indices, coming from \cite{bj1}, \cite{bj2}. The other way is by using the computation in the end of section 3, together with the fact that a free product of arbitrary Fuss-Catalan planar algebras is a Fuss-Catalan planar algebra (\cite{bj5}).

\section{Disconnected graphs}

In order to include the disjoint sum operation in the framework of previous section, we work here with slightly more general graphs. 

\begin{defi}
A precolored graph $X=(V,E)$ is a finite set $V=X_n$ together with a set $E=\{E_1,\ldots,E_p\}$ of disjoint subsets of $V\times V-\Delta_V$.
\end{defi}

Any colored graph can be regarded as a precolored graph. Conversely, for a precolored graph $X=(V,E)$ with $E=\{E_1,\ldots,E_p\}$, let $E_{p+1}$
be the set of of missing edges.
$$V\times V-\Delta_V=E_1\sqcup\ldots\sqcup E_p\sqcup E_{p+1}$$

This partition is denoted $\overline{E}$, and defines a colored graph $\overline{X} =(V, \overline{E})$. We define the Laplacian of $X$ to be that of $\overline{X}$, and the Hopf algebra associated to $X$ to be that associated to $\overline{X}$.
$$A(X)=A(\overline{X})$$

The construction of free products makes sense for precolored graphs, with the same definition, and the resulting objects are precolored graphs. However, we cannot expect Theorem 6.1 to hold in general. For instance this not the case when $X=X_n^\emptyset$ is the graph without edge. On the other hand, for $Y = X_n^\emptyset$ the free product $X*Y$ is the disjoint union of $n$ copies of $X$.
$$X*Y=X\sqcup\ldots\sqcup X$$

When $X$ is connected a free wreath product decomposition of type $A(X*Y)=A(X)*_wA(Y)$ is found in \cite{bi2}. We extend now this result to algebras $A(X)$ defined as in this paper.

\begin{defi}
A precolored graph $X=(V,E)$ with $E=\{ E_1,\ldots ,E_p\}$ is called connected if 
the oriented graph $Y=(V,F)$ given by
$$F=\bigcup_{i=1}^pE_i$$
is a connected in the usual sense. 
\end{defi}

In this definition connectedness of a usual oriented graph $Y$ means that any two vertices $i,j$ can be joined by a path of oriented edges, regardless of edge orientations.

We use the following simple fact regarding connectedness.

\begin{lemm}
If $X=(V,E)$ with $V=X_m$ is a connected oriented graph, the $*$-algebra generated by its Laplacian $d\in M_m(\c)$ contains a matrix having nonzero entries.
\end{lemm}

\begin{proof}
In order to avoid double indices we use the following notations.
$$d[1](i,j)=d_{ij}$$
$$d[*](i,j)=d^*_{ij}$$

For any $i,j$ consider a sequence of vertices realising a path from $i$ to $j$.
$$i=k_0,k_1,\ldots ,k_{n},k_{n+1}=j$$

For any $t$ we have $(k_t,k_{t+1})\in F$ or $(k_{t+1},k_t)\in F$. This means that we have $d[1](k_t,k_{t+1})=1$ or $d[1](k_{t+1},k_t)=1$, so the $(k_t,k_{t+1})$ entry of $d[1]$ or of $d[*]$ is equal to $1$.
$$d[e_t](k_t,k_{t+1})=1$$

Here each exponent $e_t$ is either $1$ or $*$. We multiply now all these equalities.
$$d[e_0](k_0,k_1)\ldots d[e_{n}](k_n,k_{n+1})=1$$

This term contributes to the $(ij)$ entry of $d[e_0]\ldots d[e_{n}]$, so we get the following inequality.
$$(d[e_0]\ldots d[e_{n}])(i,j)\geq 1$$

Consider now a sum of such matrices $d[e_0]\ldots d[e_{n}]$, one for each choice of $i,j$. This is a matrix in the $*$-algebra generated by $d$, having all entries bigger than $1$.
\end{proof}

\begin{prop}
If $X$ is a connected precolored graph we have the equality
$$A(\overline{X}*Y)=A(X*Y)$$
for any precolored graph $Y$.
\end{prop}

\begin{proof}
We use the notations $X=(X_m,E)$ and $Y=(X_n,F)$ with $E=\{E_1,\ldots , E_p\}$ and  $F= \{F_1 \ldots , F_q\}$, as in Definition 5.4. For any $r$ the Laplacian of the oriented graph $(X_m \times X_n, E_r^\circ)$ is $d_r \otimes 1_n$, where $d_r$ is the Laplacian
of $(X_m,E_r)$ and $1_n \in M_n(\c)$ is the identity matrix. Also for any $s$ the Laplacian of the oriented graph $(X_n \times X_n, F_s^\circ)$ is $\i_m \otimes \delta_s$, where $\mathbb I_m \in M_m(\c)$ is the matrix having $1$ as entries and $\delta_s$ is the Laplacian of $(X_n,F_s)$.

Consider the universal magic biunitary matrix $w\in M_{mn}(A(X_{mn}))$. Then $w$ commutes with the Laplacian of $X*Y$ if and only if it commutes
with the following matrices.
$$d_r \otimes 1_n \ , \ \mathbb I_m \otimes \delta_s, \
r=1, \ldots , p, \ s=1 \ldots q \quad (\star)$$

Also $w$ commutes with the Laplacian of $\overline{X}*Y$ if only if 
it commutes with matrices $(\star)$ and with the following matrix.
$$\left(\i_m - \sum_{i=1}^rd_i\right) \otimes 1_n$$

So assume that $w$ commutes with the family ($\star$) of matrices,
and put $d = \sum_{i=1}^r d_i$. Then $w$ commutes with $d\otimes 1_n$, and being
unitary, it also commutes with $d^*\otimes 1_n$.
By the previous lemma, there exists a non-commutative
polynomial $P(x,y)$ such that $a=P(d,d^*) \in M_m(\c)$
is a matrix having nonzero entries.  Consider now the following matrix.
$$P(d\otimes 1_n,d^*\otimes 1_n)=P(d,d^*)\otimes 1_n=a\otimes 1_n$$

Then $w$ commutes with $a \otimes 1_n$. We can regard the matrix $a \otimes 1_n$ as being the 
Laplacian of a colored graph.
By identifying all nonzero colors in this graph, we the following implication.
$$[w,a\otimes 1_n]=0\Rightarrow [w,\i_m\otimes 1_n]=0$$

Therefore $w$ commutes with $(\i_m-d)\otimes 1_n$, and we get the result.
\end{proof}

\begin{theo}
If $X$ is a connected precolored graph we have an isomorphism 
$$A(X*Y)=A(X)*_w A(Y)$$
for any precolored graph $Y$.
\end{theo}

\begin{proof}
In proof of Theorem 6.1 we have used the fact that $X$ is a colored graph (in the proof of Lemma 6.1) but we have never used the fact that $Y$ is colored. So we have $A(\overline{X}*Y)=A(X)*_wA(Y)$, and therefore Proposition 7.1 applies and concludes the proof.
\end{proof}

With a little more work one can verify that the general formula $A(X*Y)=A(X)*_w A(Y)$ holds as well for algebras $A(X)$ constructed as in \cite{bi1},
so we have a generalisation of Theorem 4.2 in \cite{bi1}. On the other hand, the algebra $A(X)$ from \cite{bi1} is expected to be of form $A(\tilde{X})$, with $A$ constructed in the spirit of this paper, and with $\tilde{X}$ being a higher combinatorial structure associated to $X$, having same symmetry group as $X$. Thus we have evidence for more general decomposition results of type $A(X*Y)=A(X)*_w A(Y)$.

As explained in section 3, useful here would be such a formula with $X,Y$ subalgebras of spin planar algebras, because this would imply Conjecture 3.1.

\begin{coro}
If $Z$ is an homogeneous precolored graph we have the equality
$$A(Z)=A(X)*_wA(X_n)$$
where $X$ is a connected component of $Z$, and $n$ is the number of components.
\end{coro}

\begin{proof}
We have $Z=X\sqcup\ldots\sqcup X$, and Theorem 7.1 applies.
\end{proof}

Together with Conjecture 3.1 and with Voiculescu's $S$-transform technique, this result reduces computation of $\mu (A(Z))$ for arbitrary homogeneous graphs to that of $\mu (A(Z))$ for connected homogeneous graphs. Indeed, we have the following computation.
\begin{eqnarray*}
\mu(A(Z))
&=&\mu (A(X)*_wA(X_n))\cr
&=&\mu (A(X))\boxtimes\mu(A(X_n))\cr
&=&\mu(A(X))\boxtimes\eta_n
\end{eqnarray*}

Here all equalities are true, except maybe for the second one, which follows from Theorem 4.1 for $n=2$, and from Conjecture 3.1 for $n\geq 3$.

It is probably useful here to compute the $S$-transform of $\eta_n$.

\begin{prop}
The $S$-transform of $\eta_n$ is given by
$$S_2(z)=\frac{1+z}{1+2z}$$
$$S_3(z)=\frac{2+2z}{1+4z+\sqrt{1+4z^2}}$$
$$S_{4+}(z)=\frac{1}{1+z}$$
where $4+$ denotes any number $n\geq 4$.
\end{prop}

\begin{proof}
The first formula follows the equality $A(X_2)=\c (S_2)$ and from Proposition 4.1. For $A(X_3)=\c (S_3)$ we use the Poincar\'e series in proof of Theorem 2.2. Substacting $1$ gives $\psi$.
\begin{eqnarray*}
\psi (z)
&=&\frac{1}{6}\left(2+\frac{3}{1-z}+\frac{1}{1-3z}\right) -1\cr
&=&\frac{1}{6}\left(\left(\frac{3}{1-z}-4\right)+\frac{1}{1-3z}\right)\cr
&=&\frac{1}{6}\left(\frac{4z-1}{1-z}+\frac{1}{1-3z}\right)\cr
&=&\frac{1}{6}\cdot\frac{4z-1-12z^2+3z+1-z}{(1-z)(1-3z)}\cr
&=&\frac{z-2z^2}{1-4z+3z^2}
\end{eqnarray*}

Thus $\chi =\chi (z)$ satisfies the following equation.
$$z=\frac{\chi-2\chi^2}{1-4\chi +3\chi^2}$$
$$z-4z\chi+3z\chi^2=\chi -2\chi^2$$
$$(2+3z)\chi^2-(1+4z)\chi +z=0$$

Together with $\chi (0)=0$, this equation gives the formula of $\chi$.
$$\Delta =1+16z^2+8z-8z-12z^2$$
$$\chi (z)=\frac{1+4z-\sqrt{1+4z^2}}{2(2+3z)}$$

This gives the formula in the statement.
\begin{eqnarray*}
S(z)
&=&\frac{1+z}{z}\cdot \frac{1+4z-\sqrt{1+4z^2}}{2(2+3z)}\cr
&=&\frac{1+z}{4z+6z^2}\cdot\frac{1+16z^2+8z-1-4z^2}{1+4z+\sqrt{1+4z^2}}\cr
&=&\frac{1+z}{4z+6z^2}\cdot\frac{8z+12z^2}{1+4z+\sqrt{1+4z^2}}\cr
&=&\frac{2+2z}{1+4z+\sqrt{1+4z^2}}
\end{eqnarray*}

The last formula follows from proof of Theorem 5.1.
\end{proof}

The formula of $S_3$ looks quite complicated, when compared to those of $S_2$ and $S_{4+}$. It is tempting to get rid of the square root, with the following change of variables.

\begin{prop}
The $S$-transform of $\eta_n$ is given by
$$S_2\left(\frac{q}{1-q^2}\right) =\frac{1+q-q^2}{1+2q-q^2}$$
$$S_3\left(\frac{q}{1-q^2}\right) =\frac{1+q-q^2}{1+2q}$$
$$S_{4+}\left(\frac{q}{1-q^2}\right) =\frac{1-q^2}{1+q-q^2}$$
where $4+$ denotes any number $n\geq 4$.
\end{prop}

\begin{proof}
The first and third formula are obtained as follows.
$$S_2\left(\frac{q}{1-q^2}\right)=\frac{1+\frac{q}{1-q^2}}{1+2\,\frac{q}{1-q^2}}=\frac{1+q-q^2}{1+2q-q^2}$$
$$S_{4+}\left(\frac{q}{1-q^2}\right) =\frac{1}{1+\frac{q}{1-q^2}}=\frac{1-q^2}{1+q-q^2}$$

For the second one, we compute first the square root.
$$\sqrt{1+4\left( \frac{q}{1-q^2}\right)^2}=\frac{1}{1-q^2}\sqrt{1+q^4-2q^2+4q^2}=\frac{1+q^2}{1-q^2}$$

This gives the following formula.
$$S_3\left(\frac{q}{1-q^2}\right) =\frac{2+2\,\frac{q}{1-q^2}}{1+4\,\frac{q}{1-q^2}+\frac{1+q^2}{1-q^2}}=\frac{2-2q^2+2q}{1-q^2+4q+1+q^2}$$

After simplification we get the formula in the statement.
\end{proof}

\section{Small graphs}

This is a technical section, using terminology and results from \cite{sms}, \cite{gr}. We will be concerned with unoriented colored graphs $X=(V,E)$. Such a graph corresponds to an oriented colored graph, when replacing each edge $(ij)$ by the pair of oriented edges $\vec{ij}$ and $\vec{ji}$.

The classification of homogeneous unoriented colored graphs with small number of vertices was started in \cite{sms}, with a complete result for graphs having $n\leq 7$ vertices. The next step was to investigate the case $n=8$. A complete result is obtained in \cite{gr} in the bicolored case, with the remark that in the arbitrary case there is just one graph left. This is the graph formed by two rectangles, which can be investigated by using techniques in this paper.

In other words, we have now complete results for all homogeneous unoriented colored graphs with $n\leq 8$ vertices. We provide here the list of all measures appearing from such graphs.

First is the dihedral series of graphs, corresponding to the case $A(X)=\c (D_m)$. This is known to contain all $n$-gons with $n\neq 4$, a certain tricolored octogon corresponding to the missing group $D_4$, plus some other graphs, like the $8$-spoke wheel. See \cite{sms}, \cite{gr}.

\begin{prop}
The spectral measure of $\c (D_m)$ is given by
$$d_{2k}=\frac{1}{4k}\left( (3k-1)\delta_0+k\delta_2+\delta_{2k}\right)$$
$$d_{2k+1}=\frac{1}{4k+2}\left( 2k\delta_0+(2k+1)\delta_1+\delta_{2k+1}\right)$$
where $k$ is such that $m=2k$ or $m=2k+1$.
\end{prop}

\begin{proof}
This follows from Proposition 2.1, see \cite{gr} for details.
\end{proof}

The other series is the Fuss-Catalan one. This contains many graphs, see \cite{sms}, \cite{gr}. The spectral measures here appear as free multiplicative convolutions between Temperley-Lieb measures in Theorem 2.2, and can be written in the following way.
$$\eta_{abc}=\eta_2^{\boxtimes\, a}\boxtimes\eta_3^{\boxtimes\, b}\boxtimes\eta_{4+}^{\boxtimes\, c}$$

We don't know how to compute this measure in the general case. We just have the following formula for the corresponding $S$-transform.

\begin{prop}
The $S$-transform of the Fuss-Catalan measure $\eta_{abc}$ is given by
$$S_{abc}\left(\frac{q}{1-q^2}\right) =\frac{(1-q^2)^c}{(1+2q)^b}\cdot\frac{(1+q-q^2)^d}{(1+2q-q^2)^a}$$
where $d$ is such that $a+b=c+d$.
\end{prop}

\begin{proof}
We use Proposition 7.3, and the fact that $\log S$ linearises $\boxtimes$.
\begin{eqnarray*}
S_{abc}\left(\frac{q}{1-q^2}\right) 
&=& S_2\left(\frac{q}{1-q^2}\right)^aS_3\left(\frac{q}{1-q^2}\right)^bS_{4+}\left(\frac{q}{1-q^2}\right)^c\cr
&=&\left(\frac{1+q-q^2}{1+2q-q^2}\right)^a
\left(\frac{1+q-q^2}{1+2q}\right)^b\left(\frac{1-q^2}{1+q-q^2}\right)^c
\end{eqnarray*}

After rearranging terms, we get the formula in the statement.
\end{proof}

After removing all dihedral and Fuss-Catalan graphs, there are just three graphs left. These are the cube, the rectangle, and the two rectangles. See \cite{sms}, \cite{gr}.

\begin{prop}
The spectral measure of the cube is given by
$$d\mu (x)=\frac{1}{8\pi}\sqrt{8x-x^2}\,dx$$
on $[0,8]$, and $d\mu (x)=0$ elsewhere.
\end{prop}

\begin{proof}
The Poincar\'e series of the cube is computed in \cite{gr}. The series is expressed there as a sum, which is summed as follows.
\begin{eqnarray*}
f(z)
&=&1+\sum_{k=1}^\infty\frac{2^{k-1}}{k+1}\,\begin{pmatrix}2k\cr k\end{pmatrix}z^k\cr
&=&1+\frac{1}{2}\sum_{k=1}^\infty\frac{1}{k+1}\,\begin{pmatrix}2k\cr k\end{pmatrix}(2z)^k\cr
&=&\frac{1}{2}\left( 1+\sum_{k=0}^\infty\frac{1}{k+1}\,\begin{pmatrix}2k\cr k\end{pmatrix}(2z)^k\right)\cr
&=&\frac{1}{2}\left( 1+\frac{1-\sqrt{1-8z}}{4z}\right)\cr
&=&\frac{1+4z-\sqrt{1-8z}}{8z}
\end{eqnarray*}

In this computation the fourth equality is easy to check, and corresponds to a well-known property of the Poincar\'e series of $TL(4+)$. We can compute now $G(\xi )=\xi^{-1}f(\xi^{-1})$.
$$G(\xi )=\frac{\xi+4-\sqrt{\xi^2-8\xi}}{8}$$

The Stieltjes formula applies, and gives the measure in the statement.
\end{proof}

\begin{prop}
The spectral measure of a rectangle is given by
$$\mu =\frac{1}{4}(3\delta_0+\delta_4)$$
and the associated universal Hopf algebra is $A=\c (\z_2\times\z_2)$.
\end{prop}

\begin{proof}
The second assertion is from \cite{sms}, and the first one follows from Proposition 4.1.
\end{proof}

\begin{prop}
The spectral measure of the graph formed by two rectangles is given by
$$d\mu (x)=\frac{1}{2}\,\delta_0+\frac{1}{\pi}\cdot \frac{\sqrt{-4+8x-x^2}}{8x-x^2}\, dx$$
on $[4-\sqrt{12},4+\sqrt{12}]$, and $d\mu (x)=0$ elsewhere.
\end{prop}

\begin{proof}
The graph $R\sqcup R$ formed by two rectangles can be investigated by using Corollary 7.1, Theorem 4.1 and Proposition 8.4. We get the following spectral measure.
\begin{eqnarray*}
\mu
&=&\mu (A(R\sqcup R))\cr
&=&\mu (A(R)*_wA(X_2))\cr
&=&\mu (\c(\z_2\times\z_2)*_w\mu(\c (\z_2)))\cr
&=&\mu (\c(\z_2\times\z_2))\boxtimes\mu(\c (\z_2))
\end{eqnarray*}

Now both $\z_2\times\z_2$ and $\z_2$ being groups acting on themselves, the corresponding $S$-transforms are given by the formula in Proposition 4.1. This gives the $S$-transform of $\mu$.
$$S(z)=\frac{1+z}{1+4z}\cdot \frac{1+z}{1+2z}$$

We can compute the $\chi$-transform, then the Poincar\'e series.
$$\chi (z)=\frac{z(1+z)}{(1+4z)(1+2z)}$$
$$z=\frac{(f-1)f}{(1+4f-4)(1+2f-2)}$$
$$f^2-f=z(4f-3)(2f-1)$$
$$(8z-1)f^2-(10z-1)f+3z=0$$

Together with $f(0)=1$, this equation gives the formula of $f$.
$$\Delta=100z^2-20z+1-96z^2+12z$$
$$f(z)=\frac{10z-1-\sqrt{4z^2-8z+1}}{16z-2}$$

We compute now the Cauchy transform $G(\xi)=\xi^{-1}f(\xi^{-1})$.
\begin{eqnarray*}
G(\xi)
&=&\frac{10\xi^{-1}-1-\sqrt{4\xi^{-2}-8\xi^{-1}+1}}{\xi(16\xi^{-1}-2)}\cr
&=&\frac{10-\xi-\sqrt{4-8\xi+\xi^2}}{16\xi-2\xi^2}\cr
&=&\frac{10-\xi-\sqrt{(\xi-4)^2-12}}{2\xi (8-\xi)}
\end{eqnarray*}

This function has a pole at $\xi=0$ with residue $1/2$, and the square root is imaginary for $\xi\in [4-\sqrt{12},4+\sqrt{12}]$. The Stieltjes formula applies, and gives the following measure.
$$d\mu (x)=\frac{1}{2}\,\delta_0+\frac{1}{\pi}\cdot \frac{\sqrt{12-(x-4)^2}}{2x(8-x)}\, dx$$

This gives the formula in the statement.
\end{proof}

The next challenging graph, mentioned in \cite{gr}, is the discrete torus at $n=9$.

\end{document}